\documentclass[12pt]{article}

\usepackage{fancyhdr, amsmath, amssymb, amsfonts, url, subfigure, dsfont, mathrsfs, graphicx, epsfig, subfigure, amsthm, tikz}

\usepackage{todonotes}

    \usepackage{epsfig}
    \usepackage{graphicx}
    \usepackage{color}
    \usepackage{wasysym}
    \usepackage[utf8]{inputenc}
    
    \usepackage[singlelinecheck=false]{caption}

	\newenvironment{blue}{\color{blue}} {\color{black}}

\newlength{\guillotine}
\newtheorem{thm}{Theorem}[section]

\newtheorem{lemma}[thm]{Lemma}

\newtheorem{prop}[thm]{Proposition}

\newtheorem{definition}[thm]{Definition}

\theoremstyle{remark}
\newtheorem{rem}[thm]{Remark}

\theoremstyle{plain}

\newtheorem*{thm*}{Theorem}

\newtheorem*{lem*}{Lemma}

\theoremstyle{definition}

\newcommand{\hp}{\hphantom}
\newcommand{\vp}{\vphantom}
\newcommand{\ph}{\phantom}

\newcommand{\NN}{\mathbb N}

\newcommand{\de}{\delta}
\newcommand{\eps}{\varepsilon}

\newcommand{\la}{\lambda}
\newcommand{\sig}{\sigma}
\newcommand{\om}{\omega}

\newcommand{\calb}{\mathcal B}
\newcommand{\calc}{\mathcal C}

\newcommand{\cali}{\mathcal I}

\newcommand{\id}{\mathrm{d}} 

\newlength{\IIlength}
\newcommand{\skewmatrix}[2]{\begin{pmatrix} 0 & #1 \\ #2 & 0 \end{pmatrix}}

\DeclareMathOperator{\Span}{Span}
\newcommand{\mobius}[2]{{\frac{#1 + #2}{1 + \overline{#2}#1}}}
	
\DeclareMathOperator{\sign}{sign}

\begin{document}

\title{Explicit examples of resonances for Anosov maps of the torus}
\author{Mark Pollicott\thanks{Supported by ERC grant 833802-resonances and EPSRC grant EP/T001674/1.} {} and Benedict Sewell\thanks{Supported by the Alfr\'ed R\'enyi Institute Young Researcher Fund.}}

\maketitle

\abstract{In \cite{BJS - anosov} Slipantschuk, Bandtlow and Just gave concrete examples of Anosov diffeomorphisms of $\mathbb T^2$ for which their resonances could be completely described.  Their approach was based on composition operators acting on analytic  anisotropic Hilbert  spaces.  
In this note we present a construction of alternative anisotropic Hilbert spaces which helps to simplify parts of their analysis and gives scope for constructing further examples.}

\section{Introduction}

In the study of chaotic diffeomorphisms, a natural class of examples are Anosov diffeomorphisms.  In fact, it is the principle of the Cohen-Gallavotti chaotic hypothesis that chaotic behaviour can be understood through the dynamics of Anosov systems \cite{cohen-gallavotti}. 

The study of Anosov dynamics is advanced by understanding various dynamical quantities, including  the resonances.
Given a map $T$, its resonances comprise a sequence (finite or converging to zero) of distinct complex numbers $(\rho_n)_{n=1}^\infty$, which give all possible exponential decay rates for the correlation function
	$$
			\int f \circ T^m g \;d\mu
		-
			\int f\, d\mu
			\int g\, d\mu, \qquad m \geq 0,
	$$
for all (sufficiently smooth) observables $f$ and $g$, and where $\mu$ is the SRB measure  (see, e.g., \cite{young} for an account of the SRB measure, and \cite[Theorem 7.11]{baladi book} for a precise statement).

Until recently, the only examples for which these resonances are completely known were given by linear hyperbolic diffeomorphisms (which represent all hyperbolic diffeomorphisms of tori up to isotopy \cite{mapping class}). These examples, including the Arnol'd CAT map $B_0: \mathbb T^2 \to \mathbb T^2$ of \cite{Arnold-Avez},
	$$
			B_0: \binom ab
		\mapsto 
			\begin{pmatrix}
				2&1\\1&1
			\end{pmatrix}
			\binom ab 
		\mod 1,
	$$
have only the trivial resonances ($0$ and $1$). On the other hand,  Adam in \cite{Adam} showed that generic small perturbations of these linear diffeomorphisms yield at least one non-trivial resonance. 
In the context of pseudo-Anosov surface homeomorphisms a description for resonances of linear pseudo-Anosov maps was recently given in \cite{Gouezel and friends}.
However, of particular interest to us are the very interesting examples are given in the striking work \cite{BJS - anosov} of Slipantshuk, Bandtlow and Just. More explicitly, they provide a family of Anosov diffeomorphisms, $(B_\lambda)$, perturbing $B_0$ above, for which the resonances $(\rho_n)_n$ (with respect to real analytic functions $f$ and $g$) are infinite and explicitly known:
	$$
		\{\rho_n\}_{n=1}^\infty \subset  \{0,1\}\cup\{\lambda^n, \overline\lambda^n \,:\, n \in \mathbb N\};
	$$
where $\lambda$ is an arbitrary complex parameter with $|\lambda|<1$.%
\footnote{This inclusion will be an equality for  generic choices of $f$ and $g$, i.e., on the complement of countably many codimension one hyperplanes.}
The  resonances of an Anosov map $T$ are calculated as the eigenvalues of its composition operator, $\mathcal C_T:f\mapsto f\circ T$, or its adjoint, the transfer operator,
acting quasi-compactly on a suitable Banach space.
In particular, we can rewrite the correlation function as $\int \calc_T^m(f) g \,d\mu - \int f\,d\mu \int g \,d\mu$ for $m \geq 0$ and then deduce that, for any $\eps > 0$, there exist polynomials $\{p_n\}_{n=1}^N$ such that
$$
			\int f \circ T^m g \;d\mu
		-
			\int f\, d\mu
			\int g\, d\mu = \sum_{n=1}^N p_n(m) \rho_n^{m} + \mathcal O \left(\eps^{m} \right), \qquad m \geq 0
	$$
	(where the degree of $p_n$ is determined by the multiplicity of $\rho_n$).
The ambient spaces, known as anisotropic spaces, have to be tailored to the diffeomorphism $T$ and their construction is non-trivial. (A description  of the myriad anisotropic spaces seen in the literature are given an overview in \cite{Demers - gentle} and a more thorough account in the survey \cite{Baladi - quest}.)

\bigskip
In this article, inspired by  \cite{BJS - anosov}, we give a new account of the resonances of $B_\lambda$ and other related examples. In particular, rather than using the spaces in \cite{BJS - anosov} (which are, in turn, based on \cite{Faure-Roy})
we introduce a new family of anisotropic Hilbert spaces
using what we call a \textit{degree function}. 
The main advantage of this construction is that  it allows us to simplify the technical analysis substantially.
Moreover, this  approach also allows us to prove new results on the resonances 
in greater generality, which we illustrate 
by  two other families, $T_\lambda$ and $T_\lambda \circ T_\mu$ in \S 2 and \S 3, respectively, where 
throughout this note, $\lambda$ and $\mu$ will denote complex parameters with $|\lambda|,|\mu|<1$.
The results on the former family appear to be new.  
The resonances of the latter family are studied empirically in an appendix of  \cite{BJS - anosov}, but  we will give a rigorous proof.

%
	
	
	
	
%

\bigskip 
We recall that a   diffeomorphism $T: \mathbb T^2 \to \mathbb T^2$ is {\it Anosov}
if there exists a continuous $DT$-invariant splitting of the tangent space $\mathcal T_\ast\mathbb T^2 = E^s\oplus E^u$ such that there exists $C>0$
and $0 < \lambda < 1$ such that $\|DT^n|E^s\| \leq C \lambda^n$ and $\|DT^{-n}|E^u\| \leq C \lambda^n$, for all $n \geq 0$.
Although the examples in this note are all Anosov, the proofs of the results are self-contained and don't depend on general properties of Anosov diffeomorphisms.

\subsection{Contents of this note}

In \S2, \S3 and \S4, respectively, we follow the general strategy of \cite{BJS - anosov} for three different families of Anosov maps $B_\lambda$, $T_\lambda$ and $T_{\lambda}\circ T_\mu$:
\begin{enumerate}
	\item[(i)] For each family, we exhibit a family of anisotropic Hilbert spaces, and show that these can be chosen to contain any pair of functions analytic on a neighbourhood of the torus.
	
	\item[(ii)] We also show that the composition operator acts compactly on these spaces (so that its spectrum gives the resonances of the map).
	
	\item[(iii)] Finally, we calculate the spectrum of this operator using a convenient, block-triangular matrix form.
	
\end{enumerate}
These results appear in the  thesis of the second author \cite{thesis}.

\section{The 
resonances of $B_\lambda$
}

The family of Anosov diffeomorphisms $B_\lambda$ 
(for $|\lambda| < 1$)
 studied in \cite{BJS - anosov} are given by so-called two-dimensional Blaschke products, originally introduced in more generality by \cite{Pujals-Shub}, where some ergodic properties were established (see also \cite{Pujals-Roeder}). More explicitly, considering
	$$
		\mathbb T^2 = \mathbb T\times \mathbb T\subset \mathbb C\times \mathbb C, \hbox{ where } \mathbb T := \{z\in\mathbb C: |z|=1\},
	$$
we have the following definition. 

\begin{definition}[$B_\lambda$] 
	Let $B_\lambda:\mathbb T^2 \to \mathbb T^2$ be given by
		$$
				B_\lambda:
					(z,w) 
			\mapsto 
				\left(
					\left(
						\frac{z+\lambda}{1 + \overline \lambda z}
					\right)
					zw,
					\left(
						\frac{z+\lambda}{1 + \overline \lambda z}
					\right)
				w
				\right).
		$$
	\end{definition}
	This family of maps analytically  perturbs the standard Arnol'd CAT map, represented 
	on $\mathbb T^2$  by
		$
			B_0: (z,w) \mapsto (z^2w,zw)
			.
		$
	%
%
The  maps 	$B_\lambda$
are Anosov and area-preserving for all $\lambda$ satisfying $|\lambda| < 1$ 
 \cite{thesis,BJS - anosov}. In particular, the resonances are well-defined and the SRB measure is just the unit area measure.

We will reprove the following result on the resonances of $B_\lambda$. This is the main result of \cite{BJS - anosov}, and we provide a new, simplified perspective.
\begin{thm}[Slipantschuk, Bandtlow and Just]
	\label{bigtheorem}
	Given $\lambda$ with $|\lambda| < 1$, there exists an  area preserving Anosov diffeomorphism for which the resonances  with respect to analytic functions $f,g: \mathbb T^2 \to \mathbb R$ take the form
		$$
			\{0,1\} \cup \{\lambda^m,\,\overline \lambda \vphantom\lambda^m\,:\,m\in \mathbb N\}.
		$$
	Moreover, each non-zero value is simple, up to coincidences in value\footnote{By this, we mean under the assumption that the $\lambda^n$ and $\overline \lambda^n$ are all distinct.}, and is otherwise semi-simple (i.e., the algebraic and geometric multiplicities coincide) of multiplicity two.
\end{thm}

The proof is based on the construction of a (non-canonical)  Hilbert space, $\mathcal H_a$, consisting of distributions on the torus, on which the composition operator $\mathcal C_{B_\lambda}: f \mapsto f \circ B_\lambda$ acts compactly and has the spectrum described in the theorem. We now describe the construction of the new  Hilbert spaces $\mathcal H_a$ we will use in the next section. 
\begin{figure}[tbh]
\centerline{
	\includegraphics[width=0.7\linewidth]{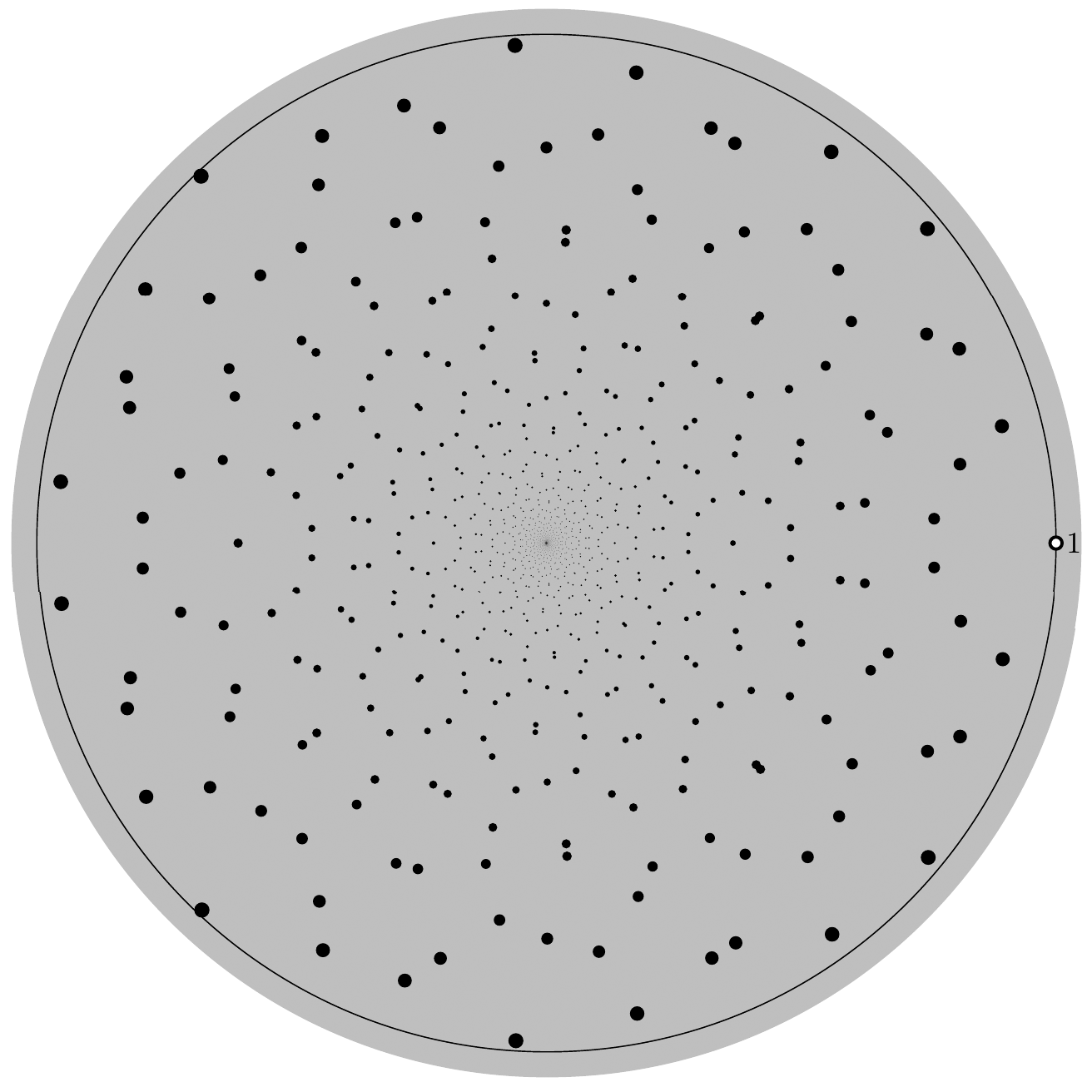}
	}
	\caption{The spectrum of $\mathcal C_{B_\lambda}$, for $\lambda = 0.99 e^{37i\pi/50}$.}
\end{figure}

\subsection{The Hilbert space $\mathcal H_a$}

All the Hilbert spaces discussed in this note are constructed using the following basic method.
Consider a complex Hilbert space $\mathcal H$ which has as an orthogonal%
\footnote{But not necessarily orthonormal.}
 basis the collection of monomials
$\{e_{m,n}\}_{(m,n)\in \mathbb Z^2}$ given by 
	$$
		e_{m,n}: (z,w) \mapsto z^mw^n.
	$$
Denoting $\langle\cdot,\cdot\rangle$ and $\|\cdot\|$ for the inner product and norm on $\mathcal H$ respectively, we have
	$$
			\left \langle
				\sum_{(m,n)\in \mathbb Z^2}
					b_{m,n}
					e_{m,n}
			,
				\sum_{(m,n)\in \mathbb Z^2}
					c_{m,n} e_{m,n}
			\right\rangle
		=
			\sum_{(m,n)\in \mathbb Z^2}
				b_{m,n}
				\overline{c_{m,n}}\,
				\|e_{m,n}\|^{2}
	$$
and
$$
\Big\| \sum_{(m,n)\in \mathbb Z^2} b_{m,n} e_{m,n} \; \Big\|^2
=
\sum_{(m,n)\in \mathbb Z^2}
|b_{m,n}|^2
\|e_{m,n}\|^2.
$$
We define $\mathcal H$ to comprise those series with finite $\|\cdot\|$ norm:
$$
		\mathcal H
	=
		\left\{
				\sum_{(m,n)\in \mathbb Z^2}
					b_{m,n} e_{m,n}
			\;\bigg|\;
				b_{m,n} \in \mathbb C,\ %
				\sum_{(m,n)\in \mathbb Z^2}
					|b_{m,n}|^2
					\|e_{m,n}\|^2
			<
				\infty
		\right\}.
$$
In particular, $\mathcal H$ is completely characterised by the values $\|e_{m,n}\|$ which we call the \textit{weights}.

\begin{rem}
	For any $a>0$, classical examples of such spaces include the Sobolev space of $a$-times weakly differentiable functions \cite[p.42]{sullivan}, which can be defined by
		$
			\|e_{m,n}\| = \big(|m| + |n| + 1\big)^{a}.
		$
	%
	Unfortunately, these spaces do not suffice for our purposes.
\end{rem}

To obtain the required properties for the composition operator acting on the Hilbert space $\mathcal H$, we need to define the weights in an anisotropic manner. In particular, taking limits along rays based at the origin, these weights decay to zero in some directions and diverge to infinity in others, and it is this behaviour which characterises the anisotropic nature of the space.

\begin{rem}
	In \cite{BJS - anosov}, after  \cite{Faure-Roy},  the authors base these weights on the eigenvectors of the  map $B_0$: i.e., for $a>0$,
		\begin{equation}
					\|e_{m,n}\|
			=  
				\exp
				\left(
					-a\left|
						\frac
							{\sqrt 5 + 1}
							2
						m + n
					\right| 
				+ 
					a\left|
						\frac
							{1 - \sqrt 5}
							2
						m + n
					\right| 					
				\right).
		\label{eqb-trionorm}
		\end{equation}
	These are a particular instance of the anisotropic spaces introduced in greater generality by Faure and Roy in \cite{Faure-Roy} and also used by Adam \cite{Adam}.
	The two essential properties of such  Hilbert spaces 
	are  that the composition operator $\mathcal C_{B_\lambda}$ acts compactly on them, and that $a>0$ can be chosen so that the space contains any given pair of functions analytic on a neighbourhood of the torus.
\end{rem}
Assuming it acts compactly, the computation of the spectrum of the composition operator acting on $\mathcal H$ 
above is to some extent  independent of the specific weights used. We therefore present  simple alternative weightings, yielding new families of anisotropic Hilbert spaces. These spaces will be particularly simple for $B_\lambda$; although  we will need a small adjustment when we consider $T_\lambda$ in the next section.

The definition of the spaces $\mathcal H_a$, appropriate to  $B_\lambda$, make use of the \textit{degree function} $\deg_1$, which we now give.
\begin{definition}[$\deg_1$, $\|\cdot\|_a$, $\mathcal H_a$]
	Let $\deg_1:\mathbb Z^2 \to \mathbb Z$ be given by
		$$
			\deg_1(m,n) := \sign(mn)\big(|m|+|n|\big),	
		$$
	where
		$$
				\hbox{\rm sign}(k) 
			= 
				\begin{cases} 
					\ph-1, & \text{if }k \geq 0;\\
					-1, & \text{if }k < 0.
				\end{cases}
		$$
 We define, for $a>0$,
		$$
				\|e_{m,n}\|_{a}
			:=
				e^{-a\deg(m,n)}.
		$$
	As described above, we let $\mathcal H_a$ be the space of series in $e_{m,n}$ with finite $\|\cdot\|_a$ norm:
	%
		$$
				\bigg\|
					\sum_{(m,n)\in \mathbb Z^2}
						b_{m,n} e_{m,n}\;
				\bigg\|^2_a
			:=
				\sum_{(m,n)\in \mathbb Z^2}
					|b_{m,n}|^2
					e^{-2a\deg_1(m,n)}.		
		$$
\end{definition}	

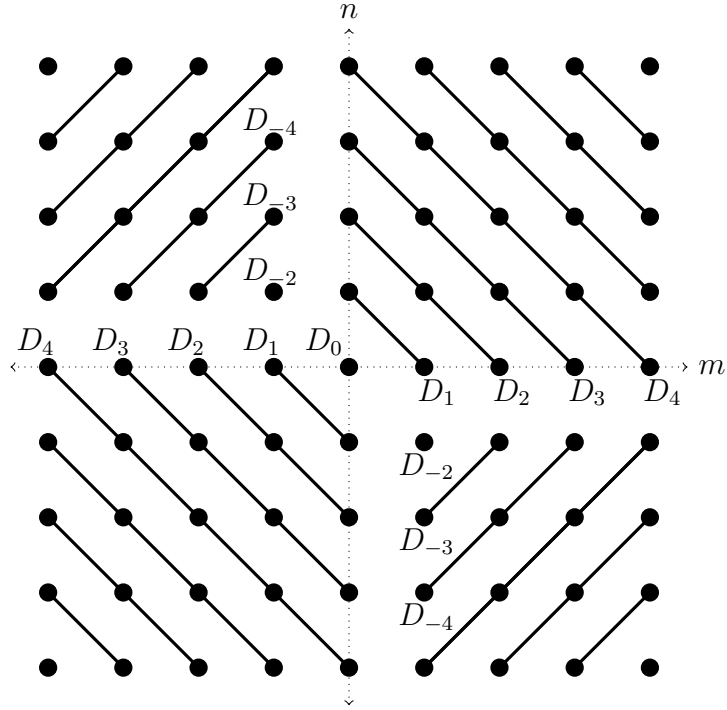
\begin{figure}[bht]
	\centering
	\begin{tikzpicture} 
	
	\draw[<->,dotted] (0,-4.5) -- (0,4.5) node [anchor = south] {$n$};
	\draw[<->,dotted] (-4.5,0) -- (4.5,0) node [anchor = west] {$m$};
	\begin{scope}[very thick]
	\foreach \x in {0,1,...,3}{
		\draw
			(\x,0) -- (0,\x);
		\draw
			(-\x,0) -- (0,-\x);
		\draw
			(\x+1,-1) -- (1,-\x-1);
		\draw
			(-\x-1,1) -- (-1,\x+1);
	}
	\foreach \x in {-4,...,4}
	\foreach \y in {-4,...,4}{
		\draw[fill=black] (\x,\y) circle[radius=0.1];
	}
	\foreach \x in {0,1,...,3}{
		\draw
			(\x,4) -- (4,\x);
		\draw
			(-\x,-4) -- (-4,-\x);
		\draw
			(\x+1,-4) -- (4,-\x-1);
		\draw
			(-\x-1,4) -- (-4,\x+1);
	}
	\end{scope}
	\foreach \x in {1,...,4}{
		\path (\x + 1/6,-1/3) node{$D_\x$};
		\path (-\x - 1/6, 1/3) node{$D_\x$};
	}
	\foreach \x in {2,...,4}{		
		\path (-1.05, \x-3/4) node{$D_{-\x}$};
		\path (1.03, -\x + 2/3) node{$D_{-\x}$};
	}
	\path (-1/3,1/3) node{$D_0$};
	\end{tikzpicture}
	\caption[The level sets of $\deg_1$]{The level sets of $\deg_1:\mathbb Z^2 \to \mathbb Z$. Here, $D_n$ denotes ${\deg_1}^{-1}(n)$.}
	\label{figb-level sets of deg 1}
\end{figure}

\noindent
	Figure \ref{figb-level sets of deg 1} shows some level sets of $\deg_1$.

The benefits of using $\mathcal H_a$ over the original family of anisotropic spaces defined by (\ref{eqb-trionorm})  can be summarized as follows.
The  proofs for compactness of the composition operators $\mathcal C: \mathcal H_a \to \mathcal H_a$ and the inclusion of analytic functions in $\mathcal H_a$ appear simpler and more direct.
Secondly, the construction permits more flexibility.  (For example,  it works also for the families $B_{\lambda,K}$ in 
the final section).
Finally, there is a clearer link between the structure of the space and the simple
 (block-triangular)
 form for the matrix of the operator.

The following result shows that any pair of analytic functions on a neighbourhood of the torus will be contained in some $\mathcal H_a$, allowing us to equate the resonances of $B_\lambda$ with the spectrum described in Theorem \ref{bigtheorem}.

\begin{prop}
	\label{prob-analytic functions in H_a}
	Let $a>0$ and suppose that $f$ is an analytic function on a neighbourhood of the poly-annulus
	$$
	P_a 
	:= 
	\left\{
	(z,w)\in \mathbb C^2
	\;\big|\;
	e^{-a}\leq|z|\leq e^a,\ %
	e^{-a}\leq|w|\leq e^a
	\right\}.
	$$
	Then $f \in \mathcal H_a$. In particular, every function analytic on a neighbourhood of $\mathbb T^2$ is contained in $\mathcal H_a$ for all sufficiently small $a$. 
\end{prop}
\begin{proof}
	Fix $a$ and let $f \in \mathcal H_a$. By construction, the Laurent series for $f$ converges absolutely on $P_a$. In particular, writing this expansion as
	\begin{equation}
	f(z,w)
	= 
	\sum_{(m,n) \in \mathbb Z^2}
	b_{m,n}z^mw^n,
	\end{equation}
	we have, by definition of $\|\cdot\|_a$,
	\begin{align*}
	\|f\|_a^2
	=
	\sum_{(m,n) \in \mathbb Z^2}
	|b_{m,n}|^2  \,
	e^{-2a\deg_1(m,n)}
	\leq
	\sum_{(m,n) \in \mathbb Z^2}
	|b_{m,n}|^2	\,
	e^{2 a(|m|+|n|)},
	\end{align*}
	which we want to show is finite. Note that
	\begin{align*}
	\sum_{m,n }
	|b_{m,n}|
	e^{a(|m|+|n|)}
	&\leq
	\sum_{m,n} 
	|b_{m,n}|
	\left(
	e^{a(m+n)}+
	e^{a(m-n)}+
	e^{a(n-m)}+
	e^{-a(m+n)}
	\right)
	\end{align*}
	is finite, since 
	(2)
	converges absolutely for all $(z,w) \in P_a$: i.e., the sums
	\begin{align*}
	&\sum_{m,n}
	|b_{m,n}|e^{a(m+n)},\qquad
	\sum_{m,n}
	|b_{m,n}|e^{a(m-n)},	
	\\&\sum_{m,n}
	|b_{m,n}|e^{a(n-m)},\qquad
	\sum_{m,n}
	|b_{m,n}|e^{-a(m+n)}							
	\end{align*}
	are each finite, since $(e^{\pm a},e^{\pm a})\in P_a$. In particular, the left hand side is square-summable, and hence $f \in \mathcal H_a$ as required.	
\end{proof}

\subsection{$\mathcal C_{B_\lambda}$ is Hilbert-Schmidt}

Since  the composition operators can be understood through their action on the  basis functions we need us estimate the corresponding 
Taylor series coefficients that appear.

\subsubsection{Estimates on Taylor coefficients}

The following definition will be used throughout.

%
\begin{definition}[$\alpha_{m,k}$]
	For all $m\in \mathbb N_0$, the following expansion converges uniformly on every disk of radius less than $|\lambda|^{-1}$:
	\begin{equation}
	\label{taylor1}
			\left(\frac{z + \lambda}{1 + \overline \lambda z}\right)^m
		=
			\sum_{k=0}^\infty
				\alpha_{m,k} z^k.
	\end{equation}

\end{definition}

\noindent
The complex coefficients $\alpha_{m,k}$ can be formulated explicitly using the Cauchy integral formula
	or Newton's identity. 
In particular,  we have $\alpha_{m,0} = \lambda^m$ for all $m \in \mathbb N_0$, and $\alpha_{0,k} = 0$ for all $k\in \mathbb N$. 

Using symmetry, one also obtains a related Taylor expansion about $\infty$ for $m \leq -1$:
	\begin{equation}
	\label{taylor2}
			\left(
				\frac{z + \lambda}{1 + \overline \lambda z}
			\right)^m
		=
			\left(
				\frac {z^{-1} + \overline\lambda}{1 + \lambda z^{-1}}
			\right)^{-m} 
		=
			\sum_{k=0}^\infty
				\overline{\alpha_{-m,k}}
				z^{-k}.
	\end{equation}
For simplicity, we adopt the notation that $\alpha_{-m,k} = \overline{\alpha_{m,k}}$ for all $m$ and $k$.

As observed in   \cite{BJS - anosov} the proof of compactness of $\mathcal C_{B_\lambda}$  reduces to estimating sums of the form
$$
\sum_{k=0}^\infty |\alpha_{m,k}|^2 e^{-2ak}
$$
for each $m \in \mathbb Z$, and for $a>0$ fixed.
In  Lemma 2.3 of  \cite{BJS - anosov} this was derived using  the Cauchy integral formula.
%
%
		%
		%
		%
%
We now present an alternative estimate, which has the advantages of being direct, simple and explicit. 

\begin{lemma}
	For all $\lambda$ and $a>0$,
		\begin{equation}
				M_{a,\lambda}
			:=
				\max_{|z| = e^{-2a}} 
					\left| 
					\frac{z+\ \lambda}{1 + \overline \lambda z}
					\right|
			<
				1.
		\label{less}
		\end{equation}
	Moreover, $M_{a,\lambda}$ satisfies, for all $m \in \mathbb Z$,
		$$
				\sum_{k=0}^\infty
					|\alpha_{m,k}|^2 e^{-2ak}
			\leq
				M_{a,\lambda}^{|m|}.
		$$
	%
\end{lemma}

\begin{proof}
	Since $|\alpha_{m,k}| = |\alpha_{-m,k}|$,  it suffices to  assume $m\geq 0$.	Since	
		$$
				\frac 1 {2\pi}
				\int_{\mathbb T^1} z^{k-j} \;|d z|
			:=
				\frac 1 {2\pi}
				\int_0^{2\pi} e^{i(k-j)\theta} d\theta
			=
			\begin{cases}
			1, & \hbox{ if } k=j;\cr
			0, & \hbox{ if } k\neq j;\cr
			\end{cases}
		$$
	we have the following, exchanging sums and integral:
		\begin{align*}
				\sum_{k=0}^\infty
					|\alpha_{m,k}|^2 e^{-2ak}	
			&=
				\frac 1 {2\pi}
				\sum_{k=0}^\infty
					\sum_{j=0}^\infty
						\alpha_{m,k}\overline{\alpha_{m,j}} e^{-2ak}
				\int_{\mathbb T^1}z^k
					z^{-j}
				\;|\id z|
			\\
			&=
				\frac 1 {2\pi}
				\sum_{k=0}^\infty
					\sum_{j=0}^\infty
						\alpha_{m,k}\overline{\alpha_{m,j}}
						\int_{\mathbb T^1}(ze^{-2a})^k
							z^{-j}
							\;|\id z|
			\\
			&=
				\frac 1 {2\pi}
				\int_{\mathbb T^1}
				\sum_{k=0}^\infty
					\alpha_{m,k}
					(ze^{-2a})^k
					\sum_{j=0}^\infty
						\overline{\alpha_{m,j}}
						z^{-j}
				\;|\id z|
			\\
			&= 
				\frac 1 {2\pi}
				\int_{\mathbb T^1}
					\left(
					\frac
					{ze^{-2a} + \lambda}
					{1 + \overline \lambda e^{-2a}z}
					\right)^m
					\left(
					\frac{z+\lambda}{1 + \overline \lambda z}
					\right)^{-m}				
					\;|\id z|,
		\end{align*}
	and a uniform estimate on this integral gives
	$$
	\begin{aligned}
	\sum_{k=0}^\infty
	|\alpha_{m,k}|^2 e^{-2ak}	
	&\leq
	\max_{|z|=1}
	\left|
	\frac
	{ze^{-2a} + \lambda}
	{1 + \overline \lambda e^{-2a}z}
	\right|^m
	{\underbrace{
		\left|
		\frac{z+\lambda}{1+\overline \lambda z} 
		\right|}_{=1}}
	^{-m}
	\\
	&=
	\max_{|ze^{2a}| = 1}
	\left|
	\frac{z+\lambda}{1+\overline \lambda z}
	\right|^m	
	=
	M_{a,\lambda}^m,
	\end{aligned}
	$$
	which proves (\ref{less}).
	Finally, elementary calculus shows that
		$$
				M_{a,\lambda}
			=
				\frac
					{|\lambda| + e^{-2a}}
					{1+e^{-2a}|\lambda|},
		$$
	leading to $M_{a,\lambda}<1$.
\end{proof}

\subsubsection{Application to $\mathcal C_{B_\lambda}$}

The previous lemma suffices to prove the following property for the composition operator $\mathcal C_{B_\lambda}$.  This immediately 
 implies compactness \cite[p.267]{conway}.
\begin{definition}[Hilbert-Schmidt, $\|\cdot\|_{\text{HS}}$]
	The Hilbert-Schmidt norm of an operator $\mathcal C:\mathcal H \to \mathcal H$ acting on a separable Hilbert space $\mathcal H$, for any orthogonal basis $\{e_i\}_{i \in \cali}$ of $\mathcal H$, is given by
		$$
				\|\mathcal C\|_{\text{\emph{HS}}}^2
			=
				\sum_{i \in \cali}
					\left(\frac
						{\|\mathcal C (e_{i})\|}
						{\|e_{i}\|}
					\right)^2.
		$$
	We say that $\mathcal C$ is \textit{Hilbert-Schmidt} if it has finite Hilbert-Schmidt norm. Note that the norm is independent of the choice of basis \cite[p.267]{conway}.
\end{definition}
%

We now show that $\mathcal C_{B_\lambda}$ has this property.

\begin{prop}
\label{compact}
	For all $a>0$, $\mathcal C_{B_\lambda}:\mathcal H_a\to \mathcal H_a$ is Hilbert-Schmidt.
\end{prop}
%

%
%

The proof of 
this proposition 
uses the following simple lemma.
\begin{lemma}\label{deg}
	For all $(m,n) \in \mathbb Z^2$, whenever $n \neq 0$,
	\begin{equation}
			\deg_{1}(m + \sign(n),n)
		\geq
			\deg_1(m,n) + 1.
	\end{equation}
	Similarly, $\deg_1(m, \sign(m)+n) \geq \deg_1(m,n) + 1$ whenever $m \neq 0$.
\end{lemma}
Although the lemma is quite intuitive (see Figure \ref{figb-level sets of deg 1})
%
%
%
we give an analytic proof for completeness.

\begin{proof}
	
	We only prove the first inequality, since the second follows by symmetry. We prove it in three cases:
	\\[4pt] \textbf{Case 1:} $mn\geq 0$. Then $(m+ \sign(n))n = mn + |n|\geq 0$ and thus
		\begin{align*}
				\deg_1(m + \sign(n),n)
			=
				\sign
				\big(\underbrace{|m+\sign(n)|}_{{|m| + 1}}+|n|\big)	
			&=
				|m| + |n| + 1
			\\
			&=
				\deg_1(m,n)+1.	
		\end{align*}		
	%
	 \textbf{Case 2:} $mn < 0$\textit{ and }$|m| > 1$. Then $mn + |n| < 0$ and thus
		\begin{align*}
				\deg_1(m + \sign(n),n)
			=
				-
				\big(\underbrace{|m+\sign(n)|}_{{|m| - 1}}+|n|\big)
			&=
				1 - |m| - |n|
			\\
			&=
				\deg_1(m,n)+1.
		\end{align*}	
	 \textbf{Case 3:} $mn < 0$\textit{ and }$|m| = 1$. The two hypotheses give $mn + |n| = 0$ and thus $\deg_1(m+ \sign(n),n) \geq 0$, whereas $mn<0$ implies that $\deg_1(m,n)\leq -1$, completing the proof.
\end{proof}
%

We now return to the proof of Proposition \ref{compact}.
\begin{proof}[Proof of Proposition \ref{compact}
]
	Fix $\lambda$ and $a>0$, and consider $\mathcal C_{B_\lambda}(e_{m,n})$. The Taylor expansions of (\ref{taylor1}) and (\ref{taylor2}) give
		\begin{align*}
			e_{m,n}\big(B_\lambda(z,w)\big)
			&=
				\left(\frac{z+\lambda}{1+\overline \lambda z}\right)^{m+n}z^mw^{m+n}
				\\
			&=
				\begin{cases}
				\rule{0pt}{20pt}\displaystyle
				\sum_{k=0}^\infty
				\alpha_{m+n,k} z^{m+\sig k}w^{m+n},	&	\text{if }m+n \neq 0;	\\
				\qquad 
				z^mw^{m+n},							&	\text{if }m+n = 0;
				\end{cases}
		\end{align*}
	where we denote $\sigma = \hbox{\rm sign}(m+n)$. That is,

	\begin{equation}
		\mathcal C_{B_\lambda}(e_{m,n})
	=
		\begin{cases}
			\displaystyle \sum_{k=0}^\infty
			\alpha_{m+n,k} \, e_{m+\sigma k,m+n},
			&\text{if } m+n \neq 0; 	
			\\
			\qquad e_{m,m+n} = e_{m,0},
			&\text{if } m+n = 0.
		\end{cases}
	\end{equation}
	Consider the case that $m+n \neq 0$. To estimate
	\begin{align}
			\left(
				\frac
					{\|\mathcal C_{B_\lambda}(e_{m,n})\|_a}
					{\|e_{m,n}\|_a}
			\right)^2		
		&=
			\sum_{k=0}^{\infty}
				|\alpha_{m+n,k}|^2
				\left(
				\frac
					{\|e_{m + \sigma k,m+n}\|_a}
					{\|e_{m,n}\|_a}
				\right)^2,
	\label{expansion}
	\end{align}
	we first bound
	\begin{equation}
			\frac
				{\|e_{m + \sigma k,m+n}\|_a}
				{\|e_{m,n}\|_a}
		=
			\exp\left[
				-a \big(
					\deg_1(m + \sigma k,m+n)
				-
					\deg_1(m,n)
				\big)
			\right],
		\label{deg11}
	\end{equation}
	for each $k \in \mathbb N_0$. To this end, we apply Lemma \ref{deg} in two different ways. Firstly, since $m+n \neq 0$, applying the lemma $k$ times gives
	$$
	\deg_1(m + \sig k,m+n) = \deg_1\big(m + \sigma k, m+n\big) \geq \deg_1(m,m+n) + k.
	$$
	Secondly, applying the lemma $|m|$ times to the right hand side gives
	$$
	\deg_1(m,m+n)
	= 
	\deg_1\big(m,|m|\hbox{\rm sign}(m)+n\big)
	\geq
	\deg_1(m,n) + |m|
	$$
	(if $m=0$, the inequality is trivial). That is,
	\begin{equation}
	\deg_1(m + \sigma  k,m+n)
	\geq
	\deg_1(m,n)
	+
	|m| + k.
	\label{eqb-deg_1 |m| + k inequality}
	\end{equation}
	Thus, by (\ref{deg11}),
	$$
	\frac
	{\|e_{m + \sigma k,m+n}\|_a}
	{\|e_{m,n}\|_a}
	\leq
	e^{-a\big(|m|+k\big)}.
	$$

	We can now bound 
	(\ref{expansion})
	 using Lemma
	\ref{less}:
	%
	\begin{align}
			\left(
				\frac
					{\|\mathcal C_{B_\lambda}(e_{m,n})\|_a}
					{\|e_{m,n}\|_a}
			\right)^2		\nonumber
		&=
			\hp{e^{-2a|m|}}
			\sum_{k=0}^\infty
				|\alpha_{m+n,k}|^2
				\left(
					\frac
						{\|e_{m + \sig k,m+n}\|_a}
						{\|e_{m,n}\|_a}
				\right)^2
		\nonumber
		\\
		&\leq
		e^{-2a|m|}
			\sum_{k=0}^\infty
				|\alpha_{m+n,k}|^2
				\;e^{-2ak}
		\nonumber
		\\
		&\leq
			e^{-2a|m|} M_{a,\lambda}^{|m+n|}
		\nonumber
		\\
		&\leq
			e^{-\de \big(|m| + |n|\big)},
			 \label{inequality2}
	\end{align}
	where $\de = \min(-\frac12\log M_{a,\lambda}, a)>0$. 
	Moreover, 
	(\ref{inequality2})
 trivially extends to the case of $m+n=0$, which is sufficient to finish the proof:
	\begin{align*}
		\|\mathcal C_{B_\lambda}\|_{\text{HS}}^2
	=
		\sum_{(m,n) \in \mathbb Z^2}
			\left(\frac
				{\|\mathcal C_{B_\lambda} e_{m,n}\|_a }
				{\|e_{m,n}\|_a}
			\right)^2
	&\leq 
		\sum_{(m,n) \in \mathbb Z^2}
			e^{-\de(|m|+|n|)}
	<
		\infty.
	\end{align*}
%
%
\end{proof}

\subsection{The spectrum of $\mathcal C_{B_\lambda}$}

As mentioned above, the calculation of the eigenvalues of $\mathcal C_{B_\lambda}$ will be independent of the weights $\|e_{m,n}\|_a$. We first give a useful definition and lemma.

\subsubsection{Block-triangular form for compact operators}

Thinking of $\mathcal C_{B_\lambda}$ as a bi-infinite matrix, we present the following definition, which generalises the notion of a block-triangular matrix, i.e., a matrix of the form
	$$
		\begin{pmatrix}
			A_{1} & 0 & 0 & \cdots & 0 \\
			\ast & A_{2} & 0 & \cdots & 0 \\
			\ast & \ast & A_{3} & \cdots & 0 \\
			\vdots & \vdots & \vdots & \ddots & \vdots \\
			\ast & \ast & \ast & \cdots & A_{n} \\
		\end{pmatrix},
	$$
where the $A_k$ are square matrices.

	This generality, although it is not required for the family $(B_\lambda)$, is convenient for when we later consider the family $(T_\lambda)$ in \S3, and is particularly so when we extend the analysis to $(T_\lambda\circ T_\mu)$ in \S4.

%
\begin{definition}[Block-triangular form]
	We say that a linear operator $\mathcal C$, acting on a Hilbert space $\mathcal H$ with orthogonal basis $\calb = \{e_{i}\}_{i \in \cali}$, has a \textit{block-triangular form} (with respect to $\calb$) if one has
		$$
				\mathcal H
			=
				\bigoplus_{k\in \mathbb Z}
				D_k
		$$
	such that, for each $k \in \mathbb Z$,
	\begin{itemize}
		\item  $D_k$ has a basis consisting of a finite (non-empty) subset of $\calb$, and
		\item  $\mathcal C(D_k) \subset \bigoplus_{j=k}^\infty D_j$.
	\end{itemize} 
	\label{lower}
\end{definition}
We now state the following result which reduces eigenvalue computations of block-triangular operators to those of their finite-dimensional blocks.
\begin{lemma}
	Suppose $\mathcal C$ and $D_k$ are as in Definition \ref{lower}, and suppose further that $\mathcal C$ is Hilbert-Schmidt. Then its non-zero eigenvalues are precisely the union of the eigenvalues for each finite rank operator $\mathcal C_k$ ($k \in \mathbb Z$):
	$$
	\mathcal C_k = \Pi_{D_k} \circ \mathcal C \circ \Pi_{D_k},
	$$
	where $\Pi_{D}$ denotes orthogonal projection onto the subspace $D$. 
	
	Moreover, if a given non-zero eigenvalue of $\mathcal C$ is an eigenvalue of only one $\mathcal C_{k}$, then its algebraic and geometric multiplicities for these two operators coincide.
	\label{lemb-eigenvalues of block triangular matrices}
\end{lemma}
This  result is quite
straightforward.   For more details see the  appendix of \cite{thesis}.

To apply this result, each of the composition operators in this note will be block-triangular with respect to $(e_{m,n})_{m,n}$, with the subspaces $D_k$ given by
\begin{equation}
D_k
=
\Span\{e_{m,n} \mid \deg_1(m,n) = k\}.
\label{eqb-D_k in terms of deg_1}
\end{equation}
Since $\deg_1(m,n) = k \implies |m|+ |n| = k$, each $D_k$ is finite dimensional, and Lemma 5 applies to any Hilbert-Schmidt operator that \textit{increases} $\deg_1$, in the following sense.

\begin{definition}[Increase]
	If $\mathcal H$ is a Hilbert space which has $(e_{m,n})_{(m,n) \in \mathbb Z^2}$ as an orthogonal basis, we say the endomorphism $\mathcal C:\mathcal H\to \mathcal H$ \textit{increases} $\deg_1$ if, for each $(m,n) \in \mathbb Z^2$, $\mathcal C (e_{m,n})$ lies in the closure of
	$$
	\Span \{e_{m',n'} \mid \deg_1(m',n') \geq \deg_1(m,n)\}, 
	$$
	i.e., $\mathcal C(D_k) \subset \bigoplus_{j = k}^\infty \,D_j$ for each $k \in \mathbb Z$, where the $D_j$ are given in
	 (\ref{eqb-D_k in terms of deg_1}).
\end{definition}

\subsubsection{Application to the spectrum of $\mathcal C_{B_\lambda}$}

We apply the above machinery to obtain the following 
useful
result, completing the proof of Theorem
\ref{bigtheorem}.

\begin{lemma}
	For all $a>0$, $\mathcal C_{B_\lambda}: \mathcal H_a \to \mathcal H_a$ has spectrum
	$$
	\{0,1\}
	\cup
	\{
	\lambda^k,\,
	\,
	\overline\lambda\vphantom\lambda^k
	\mid 
	k \in \mathbb N
	\}.
	$$
	where each non-zero eigenvalue has algebraic and geometric multiplicity equal to the frequency with which it appears in the above (in particular, they are all semi-simple).
\end{lemma}
\begin{proof}
The proof of this result is a straightforward application of Lemma \ref{lemb-eigenvalues of block triangular matrices}, recalling some details from the proof of Lemma 3.
	We first show that $\mathcal C_{B_\lambda}$ increases $\deg_1$. Considering the expansion of $\mathcal C_{B_\lambda}(e_{m,n})$ in 
	%
	%
	we have that either
	\begin{itemize}
	\item $m+n \neq 0$, and $\mathcal C_{B_\lambda}(e_{m,n})$ lies in the span of
		$
			\{ e_{m+\sig k,m+n}	\mid k \in \mathbb N_0\}
		$
	for
		$
			\sig = \sign(mn);
		$
	or 
	\item $(m,n) = (m,-m)$, and $\mathcal C_{B_\lambda}(e_{m,-m}) = e_{m,0}$.
	\end{itemize}
	Recalling (\ref{eqb-deg_1 |m| + k inequality}), in the first case we have
		\begin{equation}
			\deg_1(m + \sig k,m+n) \geq \deg_1(m,n) + |m| + k \geq \deg_1(m,n),
		\label{eqb-subsequent terms deg_1 diagonalisation}	
		\end{equation}
	and in the second case we have, from the definition,
		\begin{equation}
			\deg_1(m,0) = \deg_1(m,-m) + 3|m| \geq \deg_1(m,n).
		\label{eqb-m+n=0 case}
		\end{equation}
	Together, these show that $\mathcal C_{B_\lambda}$ increases $\deg_1$, so Lemma \ref{lemb-eigenvalues of block triangular matrices} applies.
	
	Using  the notation of that lemma, for each $j \in \mathbb Z$, the map
	$$
			\left(\mathcal C_{B_\lambda}\right)_j
		= 
			\Pi_{D_j} \circ \mathcal C_{B_\lambda} \circ \Pi_{D_j},
	$$
	can be obtained by eliminating all terms in the expansion for which the index of the basis  (i.e., $e_{m + \sig k,m+n}$) obtains a higher value of $\deg_1$ than $(m,n)$. In view of (\ref{eqb-subsequent terms deg_1 diagonalisation})--(\ref{eqb-m+n=0 case}), the only term that can remain in the $m+n \neq 0$ case is the one corresponding to $k = 0$, which remains only if $m=0$, and similarly in the $m+n = 0$ case, the single term survives only if $m=0$.

	Indeed, setting $m=0$, the zeroth term of $\mathcal C_{B_\lambda}(e_{0,n})$ is a multiple of $e_{0,n}$. More explicitly,
	$$
			\big(\mathcal C_{B_\lambda}\big)_{|n|} e_{0,n}
		=
			\alpha_{n,0}\, e_{0,n}
		=
			\begin{cases}
				\lambda^n e_{0,n},
					&\text{if }n\geq0;\\
				\overline \lambda^n e_{0,n},
					&\text{if }n<0.
		\end{cases}
	$$
	In other words, for $k < 0$, $\left(\mathcal C_{B_\lambda}\right)_k$ is the zero map, and for $k \geq 0$, it is the diagonal operator
	$$	
	\left(\mathcal C_{B_\lambda}\right)_k(e_{m,n})
	= 
	\begin{cases}
	\lambda^k e_{m,n},	& (m,n) = (0,\phantom-k);	\\
	\overline\lambda^ke_{m,n},	& (m,n) = (0,-k);	\\
	0,				& \text{otherwise.}
	\end{cases}
	$$
	Therefore, if $k> 0$, $\left(\mathcal C_{B_\lambda}\right)_k$ contributes two non-zero eigenvalues, $\lambda^k$ and $\overline\lambda^k$, and $\left(\mathcal C_{B_\lambda}\right)_0$ contributes the eigenvalue $1$.
	
	Finally, since $|\lambda|< 1$, these eigenvalues are distinct, except when $\lambda^k = \overline \lambda^k$, i.e., when $\lambda^k$ is real. In any case, since they  both appear as entries of the diagonal operator $\left(\mathcal C_{B_\lambda}\right)_k$, these eigenvalues remain semi-simple.
\end{proof}

This completes the proof of Theorem 
 \ref{bigtheorem}.

\section{The spectrum of $\mathcal C_{T_\lambda}$}

In this section, we consider a family of Anosov maps which give richer, more varied resonances.
This time, they will be perturbations of the orientation-reversing square root of the CAT map, $T_0 : \mathbb T^2 \to  \mathbb T^2$:
given by
$
T_0:(z,w) =  (zw,z)$.
\begin{definition} 
	For $\lambda$ with $|\lambda| < 1$ consider $T_\lambda: \mathbb T^2 \to \mathbb T^2$ defined by
	$$
	T_\lambda:
	(z,w) 
	\mapsto 
	\left(
	\left(
	\frac{z+\lambda}{1 + \overline \lambda z}
	\right)
	w,
	z
	\right).
	$$
\end{definition}

In this section  it is   necessary to use a slightly more complicated family of Hilbert spaces, $\mathcal H_{a,\phi}$, than in the previous section which is based on a generalisation of $\deg_1$.

The main result of this section is the following, which gives resonances for $T_\lambda$.

\begin{thm}
	For each $\lambda$ with $|\lambda| < 1$ there exists a Hilbert space $\mathcal H_{a,\phi}$ of distributions on $\mathbb T^2$, such that the composition operator $\mathcal C_{T_\lambda}: \mathcal H_{a,\phi}\to \mathcal H_{a,\phi}$ given by $\mathcal C_\lambda: f \mapsto f \circ T_\lambda$ is compact and has spectrum as follows: for $\lambda_1$ a square root of $\lambda$,
	\begin{equation}
	\{0,1\}
	\cup 	
	\{
	\omega^m \overline{\lambda_1}\vp\lambda ^n
	\mid 
	m,n\in \mathbb N_0,\,m+n\geq 1,\,\omega = \pm1
	\}.
	\end{equation}
	All non-zero eigenvalues have algebraic multiplicities as given in Lemma \ref{spectrum1}.
	Moreover, all non-zero eigenvalues are semi-simple.
	%
	\label{thmlam}
\end{thm}

This is illustrated in Figure 3 with $\lambda = 0.8 e^{31 i \pi/50}$.

\begin{figure}[thb]
\centerline{
	\includegraphics[width=0.7\linewidth]{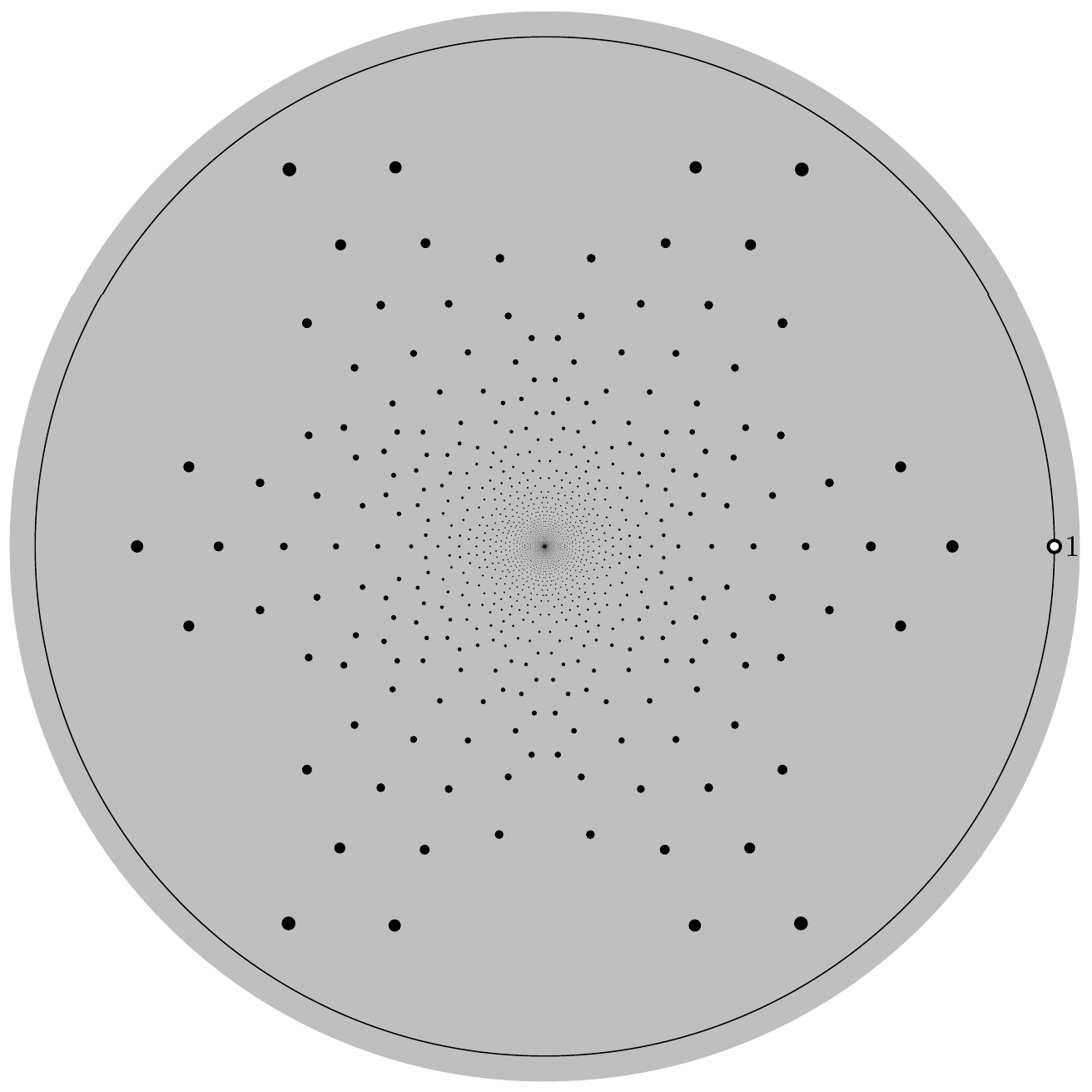}
	}
	\caption{A plot of the spectrum of $\mathcal C_{T_\lambda}$, for $\lambda = 0.8 e^{31i\pi/50}$.}
\end{figure}

\subsection{The Hilbert space $\mathcal H_{a,\phi}$}

The space $\mathcal H_{a,\phi}$ is defined analogously to $\mathcal H_a$. The weights here,  $\|e_{m,n}\|_{a,\phi},$ depend on the following simple  generalisation, $\deg_\phi$, of $\deg_1$.
\begin{definition}[$\deg_\phi$, $\|\cdot\|_{a,\phi}$, $\mathcal H_{a,\phi}$]
	For $\phi > 1$, let
		$$
				\deg_\phi(m,n) 
			:= 
				\deg_1(m,\phi^{-\sign(m,n)} n)
			=
				\begin{cases}
					\phantom{-} |m| + \phi^{-1}|n|	& \text{if }mn\geq 0;	\\
					-|m|	-\phi^{\phantom{-1}}|n| & \text{if }mn  <  0.
				\end{cases}
		$$
	For $a>0$, we write
		$$
				\|e_{m,n}\|_{a,\phi}
			:=
				e^{-a\deg_\phi(m,n)}.
		$$
	As before, this norm extends to arbitrary linear combinations of the $e_{m,n}$:
		$$
				\bigg\|
					\sum_{m,n}
						b_{m,n} e_{m,n}
				\bigg\|_{a,\phi}^2
			=
				\sum_{m,n} 
					|b_{m,n}|^2 
					e^{-2a\deg_\phi(m,n)}
					,		
		$$	
	%
\end{definition}

The following result shows that, as for  $\mathcal H_a$, the Hilbert space $\mathcal H_{a,\phi}$ can be chosen to contain 
 analytic functions on a neighbourhood of the torus.

\begin{prop}
	For $a>0$ and $\phi>1$, suppose that $f$ is an analytic function on a neighbourhood of the poly-annulus
	$$
	P_{a,\phi} := \{(z,w)\in \mathbb C^2\;|\; e^{-a}\leq|z|\leq e^a,\ e^{-a\phi}\leq|w|\leq e^{a\phi}\}.
	$$
	Then $f \in \mathcal H_{a,\phi}$. In particular, every function analytic on a neighbourhood of $\mathbb T^2$ is contained in $\mathcal H_{a,\phi}$, for all $(a,\phi)$ such that $a\phi$ is sufficiently small. 
\end{prop}
\begin{proof}
The proof 
 is very similar to that of Proposition \ref{prob-analytic functions in H_a}.
	Fix $a$, $\phi$ and $f$ as above. By construction, the expansion
	\begin{equation}
	f(z,w) = \sum_{(m,n) \in \mathbb Z^2} b_{m,n}\,z^m w^n
	\end{equation}
	converges absolutely for all $(z,w) \in P_{a,\phi}$. Also, one has the following bound from the definition of $\|f\|_{a,\phi}$, using that $-\deg_\phi(m,n) \leq |m| + \phi |n|$:
	\begin{equation}
	\|f\|_{a,\phi}^2
	:=
	\sum_{(m,n) \in \mathbb Z^2}
	|b_{m,n}|^2  \,
	e^{-2a\deg_\phi(m,n)}
	\leq
	\sum_{(m,n) \in \mathbb Z^2}
	|b_{m,n}|^2	\,
	e^{2 a(|m|+\phi|n|)}.
	\label{analytic}
	\end{equation}
	Considering the right hand side, one bounds a related sum
	\begin{align*}
	\sum_{(m,n)\in \mathbb Z^2)}
	|b_{m,n}|\,
	e^{ a(|m|+\phi|n|)}
	&\leq 
	\sum_{(m,n)\in \mathbb Z^2)}
	|b_{m,n}|				 	
	e^{a(m+\phi n)}
	+
	\sum_{(m,n)\in \mathbb Z^2)}
	|b_{m,n}|
	e^{a(m-\phi n)}
	\\
	& \ph{{}\leq{}}+
	\sum_{(m,n)\in \mathbb Z^2)}
	|b_{m,n}|
	e^{a(-m-\phi n)}
	+
	\sum_{(m,n)\in \mathbb Z^2)}
	|b_{m,n}|
	e^{a(-m+\phi n)},	
	\end{align*}
	each of which is convergent by the absolute convergence of (16)
	for all $(z,w)\in \{(e^{\pm a},e^{\pm a\phi})\}\subset P_{a,\phi}$. In particular, the sum on the left is square-summable, i.e., the sum on the right hand side of (17) is finite. Thus, $f \in \mathcal H_{a,\phi}$ as required.
\end{proof}

\subsection{$\mathcal C_{T_\lambda}$ is Hilbert-Schmidt}


To begin the proof of Theorem \ref{thmlam} we now give the following compactness result. Note that, fixing $a$ and $\lambda$, its hypothesis is satisfied for all $\phi$ sufficiently close to 1.

\begin{prop}
	\label{compact3}
	Given $\lambda$ with $|\lambda| < 1$, $a>0$ and $\phi> 1$, if
	$$
	2a(\phi-1) < -\log M_{a,\lambda},
	$$
	the composition operator $\mathcal C_{T_\lambda}:\mathcal H_{a,\phi} \to \mathcal H_{a,\phi}$ is Hilbert-Schmidt.
\end{prop}

The proof of this proposition  is similar to that of Proposition 2.10.

%
\begin{proof}
	Formally expanding
	$$
			e_{m,n}\big(T_\lambda(z,w)\big)
		=
			w^m
			\left(\mobius z w \right)^m
			z^n
	$$
	gives the following,
	 for $\sigma = \hbox{\rm sign}(m)$:
	\begin{equation}
		\mathcal C_{T_\lambda}(e_{m,n})
		=
		\begin{cases}	
			\displaystyle \sum_{k=0}^\infty
				\alpha_{m,k} \, e_{n+\sig k, m},
		& m \neq 0; 	
		\\
			\qquad e_{n,m},
		& m = 0.
		\end{cases}
	\end{equation}
	In particular, for $m \neq 0$, 
	\begin{align*}
			\left(
				\frac
					{\|\mathcal C_{T_\lambda}(e_{m,n})\|_{a,\phi}}
					{\|e_{m,n}\|_{a,\phi}} 
			\right)^2
		&=
			\sum_{k=0}^\infty
				|\alpha_{m,k}|^2
				\left(
					\frac
						{\|e_{n+\sigma k,m}\|_{a,\phi}}
						{\|e_{m,n}\|_{a,\phi}}
				\right)^2
		\\
		&=
			\sum_{k=0}^\infty
				|\alpha_{m,k}|^2
				e^{2a\big(\deg_\phi(m,n) - \deg_\phi(n + \sigma k,m)\big)}
		\\
		&=
			e^{2a\big(\deg_\phi(m,n)-\deg_\phi(n,m)\big)}
			\sum_{k=0}^\infty
				|\alpha_{m,k}|^2
				e^{2a\big(\deg_\phi(n,m) - \deg_\phi(n + \sigma  k,m)\big)}.			
	\end{align*}
	Considering first the prefactor, we find that
	\begin{align*}
			I(m,n)
		:=
			\deg_\phi(m,n) - \deg_\phi(n,m)
		&=
			\begin{cases}
		%
				\phi^{-1}(\phi - 1)\big(|m|-|n|\big),
				&\text{if }	mn \geq 0;	\\
				\phantom{\phi}^{\phantom{-1}}
				(\phi - 1)\big(|m|-|n|\big),	&\text{if }	mn  <   0.
			\end{cases}
	\end{align*}

	Also, as in the proof of Lemma 2.11, considering three cases for $\deg_\phi(n + \sig ,m) - \deg_\phi(n,m)$, we find that 
	\begin{align*}
	\deg_\phi(n + \sig,m) - \deg_\phi(n,m)
	&=
		\begin{cases}
			2|n| +  \phi^{-1}+ \phi -1,	& \text{if }mn<0\text{ and }|m| = 1;\\
			1, & \text{otherwise}.
		\end{cases}
	\\
	&\geq\rule{9.45pt}{0pt}1.
	\end{align*}

	Therefore by induction,
	$
	\deg_\phi(n + \sig k,m) - \deg_\phi(n,m) \geq k
	$
	for all $k \in \NN$.
	Thus, for all $(m,n)\in \left(\mathbb Z\backslash \{0\}\right) \times \mathbb Z$
	 (applying Lemma 2.8),
	\begin{align*}
	\left(
	\frac
	{\|\mathcal C_{T_\lambda}(e_{m,n})\|_{a,\phi}}
	{\|e_{m,n}\|_{a,\phi}} 
	\right)^2
	&\leq
	e^{2a I(m,n)}
	\sum_{k=0}^\infty
	|\alpha_{m,k}|^2 e^{-2a k }.
	\\
	&\leq
	e^{2a I(m,n)} M_{a,\lambda}^{|m|}
	\\
	&=
	\begin{cases}
	e^{2a\phi^{-1}(\phi - 1)\big(|m|-|n|\big)} 
	M^{|m|}_{a,\lambda}, 
	&\text{if }mn \geq 0;		\\
	\hp{\vp e^{\phi^{-1}}}
	e^{2a(\phi - 1)\big(|m|-|n|\big)}
	M^{|m|}_{a,\lambda},
	&\text{if }mn < 0.
	\end{cases}	
	\end{align*}
	Considering the exponents on the right hand side, if
		$$
				2a(\phi - 1)
			= 
				2a \max \big(\phi - 1, \phi^{-1}(\phi - 1)\big)
			<
				-\log(M_{a,\lambda}),
		$$
	then $\de := \min \big(2a\phi^{-1}(\phi - 1), 2a(1- \phi)-\log(M_{a,\lambda}) \big)$ is positive and satisfies
		$$
				\left(
					\frac
						{\|\mathcal C_{T_\lambda}(e_{m,n})\|_{a,\phi}}
						{\|e_{m,n}\|_{a,\phi}} 
				\right)^2
			\leq
				e^{-\de \big(|m| + |n|\big)}
		$$
	whenever $m\neq 0$. This inequality also applies in the $m = 0$ case:
	$$
	\left(
	\frac
	{\|\mathcal C_{T_\lambda}(e_{0,n})\|_{a,\phi}}
	{\|e_{0,n}\|_{a,\phi}}
	\right)^2
	=
	\left(
	\frac
	{\|e_{n,0}\|_{a,\phi}}
	{\|e_{0,n}\|_{a,\phi}}
	\right)^2
	=
	e^{-2aI(0,n)}
	\leq
	e^{-2a\phi^{-1}(\phi-1)|n|}
	\leq
	e^{-\de|n|}.
	$$
	Thus,
	\begin{equation*}
	\|\mathcal C_{T_\lambda}\|_{\text{HS}}^2
	=
	\sum_{(m,n)\in\mathbb Z^2}
	\left(
	\frac
	{\|\mathcal C_{T_\lambda}(e_{m,n})\|_{a,\phi}}
	{\|e_{m,n}\|_{a,\phi}} 
	\right)^2
	\leq
	\sum_{(m,n)\in \mathbb Z^2}
	e^{- \de (|m| + |n|)}
	<
	\infty,
	\end{equation*}
	i.e., $\mathcal C_{T_\lambda}$ is Hilbert-Schmidt, as required.
\end{proof}

	\begin{rem}
	In fact, $a(\phi-1)$ being small is necessary for $\mathcal C_{T_\lambda}$ on $\mathcal H_{a,\phi}$ to be bounded, let alone compact: For example, let $m > 0$, $n<0$. Then, considering the first term of the expansion 
	(18)
	gives
		$$
				\frac
					{\|\mathcal C_{T_\lambda}(e_{m,n})\|_{a,\phi}}
					{\|e_{m,n}\|_{a,\phi}}
			\geq
				|\lambda|^m
				\frac
					{\|e_{n,m}\|_{a,\phi}}
					{\|e_{m,n}\|_{a,\phi}}						
			=
				|\lambda|^m e^{aI(m,n)}
			=
				|\lambda|^m e^{a(\phi-1)(m + n)}.
		$$
	Thus, if $-\log|\lambda| < a(\phi-1)$, the right hand side can be made arbitrarily large.
\end{rem}

\subsection{The spectrum of $\mathcal C_{T_\lambda}$}

The following concludes the proof of Theorem
\ref{thmlam}.

\begin{lemma}
	For $\lambda$, $a$ and $\phi$ as in Proposition 3.4,
	the spectrum of $\mathcal C_{T_\lambda}:\mathcal H_{a,\phi} \to \mathcal H_{a,\phi}$ is as follows, where $\lambda_1$ is a square root of $\lambda$:
	$$
	\{0,1\} 
	\cup
	\big\{
	\om
	\lambda_1^m
	\overline\lambda_1^n
	\mid
	\om = \pm1,\ (m,n) \in \mathbb N_0^2 \setminus\{(0,0)\}
	\big\}.
	$$
	Each non-zero eigenvalue is semi-simple. Up to coincidences in value, the eigenvalues $\om\lambda_1^k$, $\om\overline\lambda_1^k$ have multiplicity
	$$
	N(k,\om)
	=
	\begin{cases}
	\left\lfloor
	\frac k2
	\right\rfloor
	+
	1,
	&	\text{if }\om = \ph-1;\\
	\left\lfloor
	\frac {k+1}2
	\right\rfloor,
	&	\text{if }\om = -1;	
	\end{cases}	
	$$
	and all other non-zero eigenvalues are simple.
	\label{spectrum1}
\end{lemma}


\begin{proof}
	The proof of this result  is analogous to the proof of Proposition 2.10.
	Recalling that $\alpha_{m,0} = \lambda^m$ for $m \in \mathbb N$, the expansion 
	 for $\mathcal C_{T_\lambda}(e_{m,n})$ reads
	$$
	\mathcal C_{T_\lambda}(e_{m,n})
	=
	\begin{cases}	
	\displaystyle \lambda^{m\hp{||}} e_{n,m} + \sum_{k=1}^\infty
	\alpha_{m,k} \, e_{n+k,m},
	&\text{if } m > 0;	
	\\
	\hp{\lambda^{|m|}}e_{n,m},
	&\text{if } m = 0;
	\\
	\rule{0pt}{20pt}
	\displaystyle
	\overline\lambda\vphantom\lambda^{|m|} e_{n,m} + \sum_{k=1}^{\infty}
	\alpha_{m,k} \, e_{n-k,m},
	&\text{if } m < 0.
	\end{cases}
	$$
	By Lemma 2.11, for any $m \neq 0$ and $k\in \mathbb N$,
	$$
	\deg_1(m,n)
	=
	\deg_1(n,m)
	<
	\deg_1\big(n + \sign(m)k,m\big).
	$$
	Since the first equality applies for $m=0$ also, this shows that $\mathcal C_{T_\lambda}$ increases $\deg_1$, and that the corresponding $\big(\mathcal C_{T_\lambda}\big)_k(e_{m,n})$ is obtained by eliminating the sums above: that is, 
	$$
	\big(\mathcal C_{T_\lambda}\big)_k 
	=
	\Pi_{D_k} \circ \mathcal C_{T_\lambda} \circ \Pi_{D_k}: e_{m,n}
	\mapsto
	\begin{cases}
	\lambda^m e_{n,m}, & \text{if }m\geq 0,\ \deg_1(m,n) = k;\\
	\overline \lambda^{|m|} e_{n,m}, & \text{if }m < 0,\ \deg_1(m,n) = k;\\
	0, & \text{otherwise;} 
	\end{cases} 	
	$$
	where $D_k := \Span\{e_{m,n}\,|\,\deg_1(m,n) = k\}$ is as before. Thus, pairing up $e_{m,n}$ and $e_{n,m}$ for $m \neq n$, one has the following block-diagonal matrix representation of $\big(\mathcal C_{T_\lambda}\big)_k$, depending on $k$:
	$$
	\big(\mathcal C_{T_\lambda}\big)_k
	\cong
	\begin{cases}
	\hphantom{\rule{0pt}{25pt}
		\displaystyle
		\bigoplus_{n=0}^{(n-2)/2}}
	\qquad\big(\,1\,\big),
	& k=0;
	\\
	\rule{0pt}{25pt}
	\displaystyle
	\bigoplus_{n=0}^{(k-2)/2}
	\skewmatrix{\lambda^n}{\lambda^{k-n}}
	\oplus 
	\skewmatrix{\overline\lambda{\vphantom\lambda}^n}{\overline\lambda{\vphantom\lambda}^{k-n}} 
	\ \oplus\;
	\big(\lambda^{k/2}\big)
	\oplus 
	\big(\overline\lambda{\vphantom\lambda}^{k/2}\big),
	& k \in 2\mathbb N;
	\\
	\rule{0pt}{25pt}
	\displaystyle
	\bigoplus_{n=1}^{(k-1)/2}
	\skewmatrix{\lambda^n}{\lambda^{k-n}}
	\oplus 
	\skewmatrix{\overline\lambda{\vphantom\lambda}^n}{\overline\lambda{\vphantom\lambda}^{k-n}},
	& k \in 2\mathbb N - 1;
	\\
	\rule{0pt}{25pt}
	\displaystyle
	\bigoplus_{\substack{n=1\\\hphantom{(k-1)/2}}}^{k-1}
	\skewmatrix{\lambda^n}{\overline\lambda{\vphantom\lambda}^{k-n}},
	& k<0.
	\end{cases} 
	$$
	Applying Lemma \ref{lemb-eigenvalues of block triangular matrices} and counting multiplicities, the non-zero eigenvalues of $\mathcal C_{T_\lambda}$ and their multiplicities are precisely those given in the statement of the lemma. In particular, each non-zero eigenvalue is semi-simple, since the $\big(\mathcal C_{T_\lambda}\big)_k$ are diagonalisable and do not share eigenvalues when $\lambda$ is non-zero.
\end{proof}

\section{The spectrum of $\mathcal C_{T_\lambda \circ T_\mu}$}
Now that we have  established the machinery for $T_\lambda$, the following result for $T_\lambda \circ T_\mu$ 
(with $|\lambda|, |\mu| < 1$) will be very easy to prove. Again, we note that this family of examples appears in an appendix of \cite{BJS - anosov}, where their resonances are announced and numerically studied. We now  provide a rigorous argument.

\begin{thm}
	For $\lambda,\mu$ with $|\lambda|, |\mu| < 1$ and $\mathcal H_{a,\phi}$ defined as above, if $a>0$ and $\phi>1$ satisfy
	\begin{equation}
	2a(\phi - 1) < -\log \big(\max( M_{a,\lambda}, M_{a,\mu})\big),
	\end{equation}
	then $\mathcal C_{T_\lambda\circ T_\mu} = \mathcal C_{T_\mu}\circ \mathcal C_{T_\lambda}$ acts compactly on $\mathcal H_{a,\phi}$ and has spectrum
	$$
	\{0,1\}
	\cup 
	\big\{ 
	\lambda^m
	\mu^n,\ %
	\lambda^m\overline\mu\vp\mu^n,\ %
	\overline \lambda^m\overline\mu\vp\mu^n,\ %
	\overline \lambda^m\mu^n
	\mid (m,n) \in \mathbb N_0^2 \setminus \{(0,0)\}
	\big\}.	
	$$
	Moreover, all non-zero eigenvalues are simple, up to coincidences in value.
\end{thm}

\begin{figure}[ht]
\centerline{
	\includegraphics[width=0.7\linewidth]{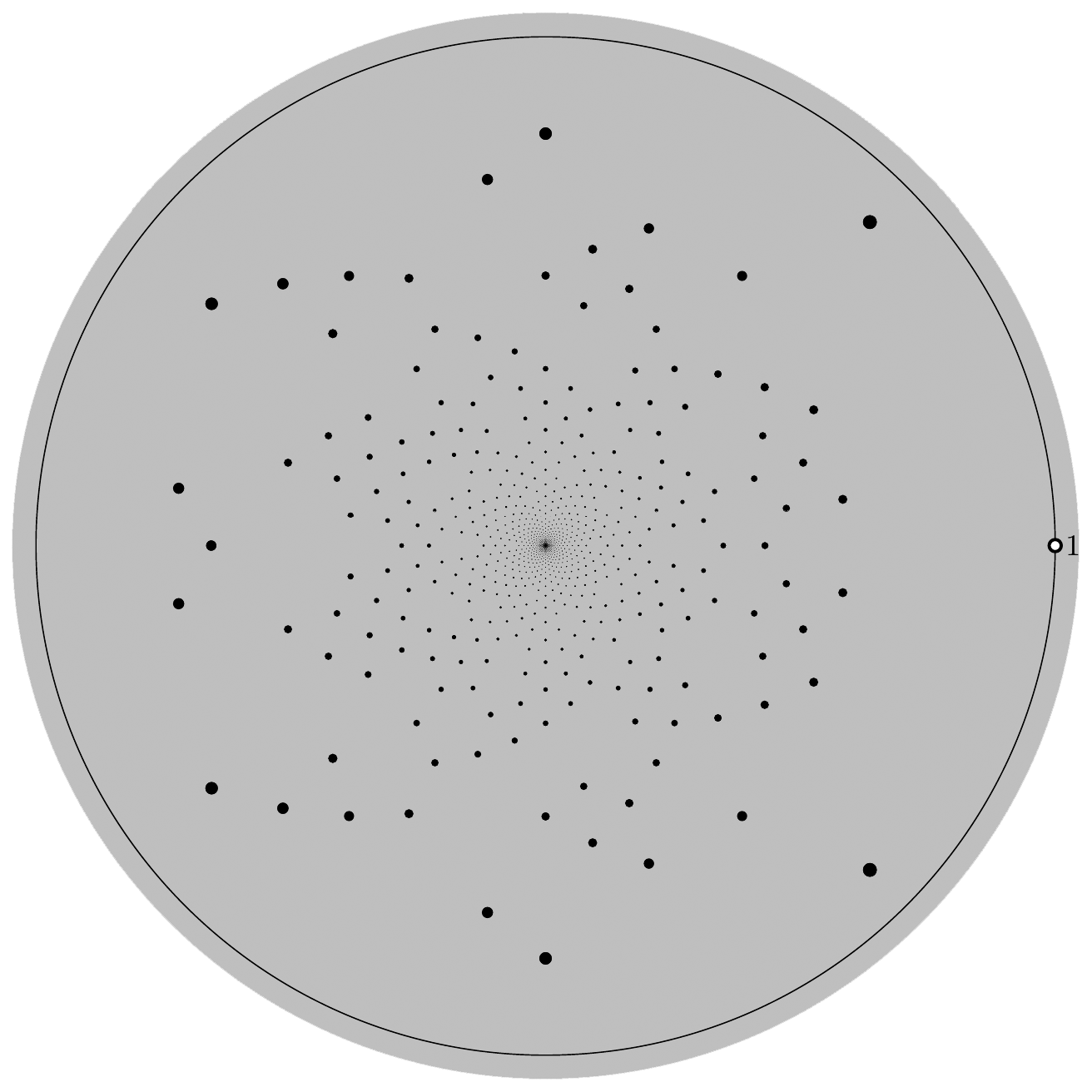}
	}
	\caption{The spectrum of $\mathcal C_{T_\lambda\circ T_\mu}$, for $\lambda = 0.9 e^{i\pi/4}$, $\mu = 0.65 e^{6 i\pi/5}$.}
\end{figure}

\subsection{$\mathcal C_{T_\lambda\circ T_\mu}$ is trace-class}
To begin the proof of Theorem 4.1
one has the following, simple corollary of Proposition 3.5.

We first recall \cite[p.267]{conway} that being trace-class is a stronger property than being Hilbert-Schmidt, and that an operator is trace-class if and only if it is the composition of two Hilbert-Schmidt operators.

\begin{lemma}
	For $\lambda,\mu,a,\phi$ as in Theorem 4.1 the operator 
	$\mathcal C_{T_\lambda\circ T_\mu}:\mathcal H_{a,\phi} \to \mathcal H_{a,\phi}$ is trace-class.
\end{lemma}
\begin{rem}
	Since $B_\lambda = T_0 \circ T_\lambda$ for all $\lambda$, this shows that $\mathcal C_{B_\lambda}$ is trace-class as an operator on $\mathcal H_{a,\phi}$.
\end{rem}
\begin{proof}
By the hypothesis (19) 
Proposition 3.5
	 applies twice to show that $\mathcal C_{T_\lambda}$ and $\mathcal C_{T_\mu}$ are each Hilbert-Schmidt on $\mathcal H_{a,\phi}$. Thus $\mathcal C_{T_\lambda\circ T_\mu} = \mathcal C_{T_\mu}\circ \mathcal C_{T_\lambda}$ is the composition of two Hilbert-Schmidt operators, hence trace-class.
\end{proof}

\subsection{The spectrum of $\mathcal C_{T_\lambda}\circ \mathcal C_{T_\mu}$}

The calculation of the spectrum likewise follows simply from that of the previous section. This uses the following lemma, which naturally extends the corresponding intuitive result for block-triangular matrices in finite dimensions:

\setlength\arraycolsep{2.5pt}
	$$
			\begin{pmatrix}
			A_{1} & 0 & 0 & \cdots & 0 \\
			\ast & A_{2} & 0 & \cdots & 0 \\
			\ast & \ast & A_{3} & \cdots & 0 \\
			\vdots & \vdots & \vdots & \ddots & \vdots \\
			\ast & \ast & \ast & \cdots & A_{n} \\
			\end{pmatrix}
			\begin{pmatrix}
			B_{1} & 0 & 0 & \cdots & 0 \\
			\ast & B_{2} & 0 & \cdots & 0 \\
			\ast & \ast & B_{3} & \cdots & 0 \\
			\vdots & \vdots & \vdots & \ddots & \vdots \\
			\ast & \ast & \ast & \cdots & B_{n} \\
			\end{pmatrix}
		=
			\begin{pmatrix}
			A_1B_{1} & 0 & 0 & \cdots & 0 \\
			\ast & A_2B_{2} & 0 & \cdots & 0 \\
			\ast & \ast &  A_3B_{3} & \cdots & 0 \\
			\vdots & \vdots & \vdots & \ddots & \vdots \\
			\ast & \ast & \ast & \cdots &  A_nB_{n} \\
			\end{pmatrix},			
	$$
where $A_k$ and $B_k$ are square matrices of the same size for each $k$.

The proof of the lemma, like that of Lemma 2.13, is a simple extension of the finite case and we omit it (see \cite{thesis}).

\begin{lemma}
	Let $\mathcal H$ be a Hilbert space such that $\{e_{m,n}\}_{(m,n) \in \mathbb Z^2}$ is an orthogonal basis, and let $\mathcal C_1, \mathcal C_2: \mathcal H \to \mathcal H$ increase $\deg_1$. Then $\mathcal C_1\circ \mathcal C_2$ increases $\deg_1$ and satisfies, for each $k$,
	\begin{equation}
	\big(\mathcal C_1\circ \mathcal C_2\big)_k
	=
	\big(\mathcal C_1\big)_k \circ \big(\mathcal C_2\big)_k.
	\end{equation}
\end{lemma}

%

We now apply this lemma to give the resonances of $T_\lambda\circ T_\mu$.

\begin{lemma}
	For each $\lambda, \mu$, 
	the spectrum of $\mathcal C_{T_\lambda\circ T_\mu}$ is given by
	\begin{eqnarray}
	\{0,1\}
	\cup 
	\big\{ 
	\lambda^m
	\mu^n,\ %
	\lambda^m\overline\mu\vp\mu^n,\ %
	\overline \lambda^m\overline\mu\vp\mu^n,\ %
	\overline \lambda^m\mu^n
	\mid (m,n) \in \mathbb N_0^2 \setminus \{(0,0)\}
	\big\}.	
	\end{eqnarray}
	Moreover, each non-zero eigenvalue has algebraic multiplicity equal to the frequency with which it appears in (21).
\end{lemma}

\begin{proof}	Applying Lemmas 4.4
	and 2.3 
	reduces the proof to a consideration of the eigenvalues of $\big(\mathcal C_{T_\lambda \circ T_\mu}\big)_k = \big(\mathcal C_{T_\lambda}\big)_k \circ \big(\mathcal C_{T_\mu}\big)_k.$ We recall from the proof of Lemma 3.7
	that, for $k = \deg_1(m,n)$,
	$$
	\big(\mathcal C_{T_\mu}\big)_k(e_{m,n}) 
	=					
	\Pi_{D_k} \circ \mathcal C_{T_\mu}(e_{m,n}) 
	=
	\begin{cases}
	\mu^{m\ph{||}} e_{n,m}, &	\text{if }m\geq 0;
	\\
	\overline\mu^{|m|} e_{n,m}, &	\text{if }m < 0.
	\end{cases}	
	$$
	Thus, for $k = \deg_1(m,n) = \deg_1(n,m)$,
	$$
	\big(\mathcal C_{T_\lambda \circ T_\mu}\big)_k(e_{m,n})
	=
	\begin{cases}
	\mu^{m\ph{||}}\lambda^{n\ph{||}}
	e_{m,n},
	&
	\text{if }m\geq 0,\ n\geq 0;
	\\
	\overline\mu^{|m|}\lambda^{n\phantom{||}}e_{m,n},	
	&
	\text{if }m < 0,\ n \geq 0;
	\\
	\mu^{m\ph{||}}\overline \lambda^{|n|}e_{m,n},
	&
	\text{if }m \geq 0,\ n < 0;
	\\
	\overline\mu^{|m|}\overline \lambda^{|n|}e_{m,n},
	&
	\text{if }m<0,\ n < 0.
	\end{cases}	
	$$
	That is, each $\big(\mathcal C_{T_\lambda \circ T_\mu}\big)_k$ is diagonal. Since the prefactor of $e_{m,n}$ is unique (up to coincidences in value), this shows that the spectrum is given by
	$$
	\{0,1\}
	\cup 
	\big\{ 
	\lambda^m
	\mu^n,\ %
	\lambda^m\overline\mu\vp\mu^n,\ %
	\overline \lambda^m\overline\mu\vp\mu^n,\ %
	\overline \lambda^m\mu^n
	\mid (m,n) \in \mathbb N_0^2 \setminus \{(0,0)\}
	\big\},
	$$
	and that the non-zero eigenvalues are simple, up to coincidences in value (e.g. if $\lambda$, $\mu$ and $\mu/\lambda$ are non-zero and have arguments which are irrational multiples of $\pi$).
\end{proof}

This completes the proof of Theorem 4.1.

\setlength\arraycolsep{4pt}

\section{Final comments}

1. 	The methods of \S 2  naturally extend to the following families of diffeomorphisms
	$B_{\lambda,K}: \mathbb T^2 \to \mathbb T^2$
	 indexed by $K\in\mathbb N$ and 
	$\lambda$ with $|\lambda|<1$:
		$$
				B_{\lambda,K}:
					(z,w) 
			\mapsto 
				\left(
					\left(
						\frac{z + \lambda}{1 + \overline \lambda z}
					\right)^{K^2+1}
					w^K,
					\left(
					\frac{z + \lambda}{1 + \overline \lambda z}
					\right)^K
					w
				\right),
		$$
	which can be considered, for each $K$, as a perturbation of the hyperbolic linear auto\-morphism given by (on $\mathbb T^2$ or $\mathbb R^2/\mathbb Z^2$ respectively)
		$$
				B_{0,K}:(z,w)\mapsto \big(z^{K^2}zw^K,\,z^Kw\big) 
			\qquad \text{or} \qquad
			\left(\begin{matrix} x \\ y \end{matrix}\right) \mapsto \begin{pmatrix}
			K^2 & \!\!K \\
			K^{\ph2}	& \!\!1
			\end{pmatrix}\left(\begin{matrix} x \\ y \end{matrix}\right) \mod 1.
		$$

	However, since the resonances of $B_{\lambda,K}$ equal those of $B_{\lambda^K}$, these families contribute nothing new to the variety of spectra presented here.

\bigskip

\noindent
2.	In \S 3 one could again extend the analysis to related families of examples: i.e., for $K\in\mathbb N$ and $|\lambda| < 1$, consider
	$$
		T_{\lambda,K}:
		(z,w) 
		\mapsto 
		\left(
		\left(
		\frac{z + \lambda}{1 + \overline \lambda z}
		\right)^K
		w,
		z
		\right),
	$$
	perturbing, for each $K$, the hyperbolic linear automorphism
	$$
	T_{0,K}:(z,w) \mapsto (z^Kw,z) 
	\qquad 
	\longleftrightarrow
	\qquad
	\left(\begin{matrix} x \\ y \end{matrix}\right) \mapsto
	\begin{pmatrix}
	K & 1 \\
	1	& 0
	\end{pmatrix}
	\left(\begin{matrix} x \\ y \end{matrix}\right)
	,
	$$
	the orientation-reversing square root of $B_{0,K}$. However again, we would find that that the spectrum of $T_{\lambda,K}$ equals that of $T_{\lambda^K}$, so these families contribute nothing extra in variety.

\bigskip
\noindent
3. 
To see why we introduced $\mathcal H_{a,\phi}$ in \S 3, we exhibit the following negative result, which shows that $\mathcal C_{T_\lambda}$ does not act compactly on either $\mathcal H_a$ or the anisotropic space used in \cite{BJS - anosov}, for any non-zero $\lambda$. 
\begin{prop}
	Suppose that $\mathcal H$ is a Hilbert space which has $\{e_{m,n}\}_{m,n}$ as an orthogonal basis, and satisfies, for all $(m,n) \in \mathbb Z^2$,
		$$
			\|e_{m,n}\| = \|e_{n,m}\|.
		$$
	Then, $\mathcal C_{T_\lambda}$ is not compact on $\mathcal H$, for any $\lambda \neq 0$.
\end{prop}
\begin{proof}
	Fix $m \in \mathbb N$ and $\lambda$. Then, recalling (18), we have
		$$
				\mathcal C_{T_\lambda}(e_{m,n})
			=
				\la^m e_{n,m}
			+
				\sum_{k=1}^\infty
					\alpha_{m,k} e_{n+k,m},
		$$
	and thus, by orthogonality, 
	\begin{equation}
			\|\mathcal C_{T_\lambda}(e_{m,n})\|^2
		\geq
			|\lambda|^{2m} \|e_{m,n}\|^2.
		\label{eqb-not compact it isnt}
	\end{equation}
	If $\mathcal C_{T_\lambda}$ is compact, it maps the sequence $(e_{m,n}/\|e_{m,n}\|)_{n=1}^\infty$, which weakly converges to zero, onto one which converges to zero in $\mathcal H$. But this contradicts (\ref{eqb-not compact it isnt}), so it is not compact.
%
\end{proof}


\begin{thebibliography}{150}
	
	\bibitem{Adam} 
	A. Adam.
	Generic non-trivial resonances for Anosov diffeomorphisms.
	\textit{Nonlinearity} \textbf{30} (2017), no. 3, 1146--64.
	\verb|doi:10.1088/1361-6544/aa59a9|
	
	
	\bibitem{Arnold-Avez} V. I. Arnold and A. Avez.
	\textit{Probl\`emes ergodiques de la mécanique classique.}
	Gauthier-Villars: Paris, 1967.
	
	
	\bibitem{baladi book}
	V. Baladi.
	\textit{Dynamical zeta functions and dynamical determinants for hyperbolic maps:
		A functional approach.} A Series of Modern Surveys in Mathematics, 3rd Series, \textbf{68}. Springer: Cham, 2018.
	ISBN:978-3-319-77660-6
	
	\bibitem{Baladi - quest}
	V. Baladi.
	The quest for the ultimate anisotropic Banach space.
	\textit{J. Stat. Phys.} \textbf{166} (2017), no. 3--4, 525--557. \verb|doi:10.1007/s10955-016-1663-0|
	
	\bibitem{bandtlow-jenkinson}
	O. Bandtlow and O. Jenkinson. On the Ruelle eigenvalue sequence.
	\textit{Ergodic Theory and Dynam. Systems} \textbf{28} (2008), no. 6, 1701--1711.
	\verb|doi:10.1017/S0143385708000059|
	
	\bibitem{BJS - expand 2}
	O. Bandtlow, W. Just and J. Slipantschuk.
	Spectral structure of transfer operators for expanding circle maps.
	\textit{Ann. Inst. H. Poincaré Anal. Non Linéaire} \textbf{34} (2017), no. 1, 31--43.
	\verb|doi:10.1016/j.anihpc.2015.08.004|
	
	
	\bibitem{bandtlow-naud}	
	O. Bandtlow and F. Naud.
	Lower bounds for the Ruelle spectrum of analytic expanding circle maps.
	\textit{Ergodic Theory Dynam. Systems} \textbf{39} (2019), no. 2, 289--310.
	\verb|doi:10.1017/etds.2017.29|
	
	
	
	\bibitem{cohen-gallavotti}
	E. G. D. Cohen and G. Gallavotti.
	Dynamical Ensembles in Nonequilibrium Statistical Mechanics.
	\textit{Phys. Rev. Lett.} \textbf{74} (1995), 2694--2697.
	\verb|doi:10.1103/PhysRevLett.74.2694|
	
	\bibitem{conway}
	J. B. Conway.
	\textit{A course in functional analysis.}
	Second edition. Graduate Texts in Math. \textbf{96}.
	Springer-Verlag: New York, 1990.
	\verb|doi:10.1007/978-1-4757-4383-8|

	\bibitem{cowen}
	C. C. Cowen.
	\textit{Composition Operators on Spaces of Analytic Functions II}.
	Lecture notes: Spring School of Functional Analysis, Rabat, 19--21 May 2009.
	\verb|math.iupui.edu/~ccowen/Talks/CompOp0905slidesB.pdf|
	
	
	\bibitem{Demers - gentle}
	M. F. Demers.
	A gentle introduction to anisotropic Banach spaces.
	\textit{Chaos Solitons Fractals} \textbf{116} (2018), 29--42.
	\verb|doi:10.1016/j.chaos.2018.08.028|
	
	\bibitem{durrett}
	R. Durrett.
	\textit{Probability: theory and examples.}
	Fifth edition. Cambridge Series in Statistical and Probabilistic Mathematics \textbf{49}. Cambridge University Press: Cambridge, 2019.
		\verb|doi:10.1017/9781108591034|
	
	\bibitem{mapping class}
	B. Farb and D. Margalit.
	\textit{A primer on mapping class groups.}
	Princeton Mathematical Series, \textbf{49}. Princeton University Press: Princeton, 2012. ISBN:978-0-691-14794-9
	
	\bibitem{Faure-Roy}
	F. Faure and N. Roy.
	Ruelle-Pollicott resonances for real analytic hyperbolic maps.
	\textit{Nonlinearity} \textbf{19} (2006), no. 6, 1233--1252.
	\verb|doi:10.1088/0951-7715/19/6/002|
	
	\bibitem{Gouezel and friends}
	F. Faure, S. Gouëzel and E. Lanneau.
	Ruelle spectrum of linear pseudo-Anosov maps.
	\textit{J. Éc. polytech. Math.} \textbf{6} (2019), 811--877.
	\verb|doi:10.5802/jep.107|
	
	
	\bibitem{naud}
	F. Naud.
	The Ruelle spectrum of generic transfer operators.
	\textit{Discrete Contin. Dynam. Systems} \textbf{32} (2012), no. 7, 2521--2531.
	\verb|doi:10.3934/dcds.2012.32.2521|
	
	

	\bibitem{Pujals-Roeder} 
	E. R. Pujals and R. K. W. Roeder.
	Two-dimensional Blaschke products: degree growth and ergodic consequences.
	\textit{Indiana Univ. Math. J.} \textbf{59} (2010), no. 1, 301--325.
	
	\bibitem{Pujals-Shub} E. R. Pujals and M. Shub.
	Dynamics of two-dimensional Blaschke products.
	\textit{Ergodic Theory Dynam. Systems} \textbf{28} (2008), no. 2, 575--85.
		\verb|doi:10.1017/S0143385707000752|
	
	
	\bibitem{rudin}
	W. Rudin.
	\textit{Functional analysis.} Second edition. International Series in Pure and Applied Mathematics. McGraw-Hill: New York, 1991. ISBN:0-07-054236-8
	
%
	
	
	\bibitem{thesis}
		B. Sewell.
		\textit{Equidistribution of infinite interval substitution schemes, explicit resonances of  Anosov  toral  maps,  and  the  Hausdorff  dimension  of  the  Rauzy  gasket.} PhD thesis. Available from \verb|www.warwick.ac.uk/fac/sci/maths/people/staff/sewell|
	
	\bibitem{shapiro}
	J. H. Shapiro.
	\textit{Composition operators and classical function theory.}
	Universitext: Tracts in Mathematics. Springer-Verlag: New York, 1993. ISBN:0-387-94067-7
	
	
	
	\bibitem{BJS - expand 1}
	J. Slipantschuk, O. Bandtlow and W. Just.
	Analytic expanding circle maps with explicit spectra.
	\textit{Nonlinearity} \textbf{26} (2013), no. 12, 3231--45.  \verb|doi:10.1088/0951-7715/26/12/3231|
	
	\bibitem{BJS - anosov}
	J. Slipantschuk, O. Bandtlow and W. Just.
	Complete spectral data for analytic Anosov maps of the torus.
	\textit{Nonlinearity} \textbf{30} (2017), no. 7, 2667--86. \verb|doi:10.1088/1361-6544/aa700f|
	
	\bibitem{sullivan}
	T. J. Sullivan.
	\textit{Introduction to uncertainty quantification.}
	Texts in Applied Mathematics \textbf{63}. Springer: Cham, 2015.
	\verb|doi:10.1007/978-3-319-23395-6|

	
	

	\bibitem{young}
	L.-S. Young.
	What are SRB measures, and which dynamical systems have them?
	\textit{J. Stat. Phys.} \textbf{108} (2002), no. 5-6, 733--754.
	\verb|doi:10.1023/A:1019762724717|
	
\end{thebibliography}
\end{document}